\documentclass[11pt]{article}%{report}
\usepackage{amssymb} \usepackage{amsmath} \usepackage{makeidx} \usepackage{esint}
\usepackage[dvips]{graphicx, color}

\newcounter{counterpar}
\newcounter{ctr}
\newtheorem{definition}[ctr]{Definition}
\newtheorem{theorem}[ctr]{Theorem}
\newtheorem{proposition}[ctr]{Proposition}
\newtheorem{lemma}[ctr]{Lemma}

\newtheorem{corollary}[ctr]{Corollary}

\newcommand{\startproof}{\noindent{\bf Proof}\quad }
\newcommand{\stopproof}{\hfill$\blacksquare$\\}%{\hfill$\Box$\\}

\newcommand{\R}{{\mathbb R}}
\title{Dissipative continuous Euler flows in two and three dimensions
}
\author{A.~Choffrut, C.~De Lellis, L.~Sz\'ekelyhidi Jr.}
\date{\today}
\begin{document}
\maketitle
\begin{abstract}
In \cite{DS3}
two of these authors construct dissipative continuous (weak) solutions to the incompressible Euler equations 
on the three-dimensional torus $\mathbb T^3$.
The building blocks in their proof are Beltrami flows,
which are inherently three-dimensional.
The purpose of this note is to show that the techniques can nevertheless be adapted to the two-dimensional case.
\end{abstract}

\tableofcontents

\bigskip

\hrule

\medskip

\section{Introduction}
In this paper we will take
\[
d=2\quad {\rm or}\quad 3\,
\]
and we will consider the incompressible Euler equations
\begin{equation}
\partial_tv +{\rm div}\,(v\otimes v)+\nabla p=0,\qquad {\rm div}\,v=0\, .
\label{eq:Euler}
\end{equation}

The following Theorem was proved for $d=3$ in the paper \cite{DS3}
and we will show in this note that:
\begin{itemize}
\item the same statement holds for $d=2$ and indeed the proof of \cite{DS3}
can be suitably modified to handle this case as well;
\item in $d=3$ the approach yields nonetheless solutions which are genuinely
$3$-dimensional (cp. with Theorem \ref{theorem:3D-flows-are-not-2D}).
\end{itemize}

\begin{theorem}
\label{theorem:dissipative-continuous}
Let $T>0$.
Assume $e\colon [0,T]\rightarrow \mathbb R$ is a smooth positive function.
Then, there is a continuous vector field $v\colon \mathbb T^d\times [0,T]\rightarrow \mathbb R^d$ 
and a continuous scalar field $p\colon \mathbb T^d\times [0,T]\rightarrow \mathbb R$ 
solving the incompressible Euler equations
in the sense of distributions and such that
\begin{equation}
e(t)=\int_{\mathbb T^d}|v|^2(x,t)\,dx,\qquad t\in [0,T].
\label{eq:e(t)=int|v|2}
\end{equation}
\end{theorem}

If $(v,p)$ is a $C^1$ solution of \eqref{eq:Euler}, we can scalar multiply the first equation by $v$ and use
the chain rule to derive the identity
\[
\partial_t \frac{|v|^2}{2} + {\rm div}_x \left(\left(\frac{|v|^2}{2} + p\right) v \right) = 0\, .
\]
Integrating this last equality in space we then derive the conservation of the total kinetic energy
\begin{equation}\label{e:energy}
\frac{d}{dt}  \int_{\mathbb T^d} |v|^2 (x,t)\, dx = 0\, .
\end{equation}
Thus, classical solutions of the incompressible Euler equations are energy conservative and 
are therefore ``rigid'' compared to the continuous solutions, for which Theorem \ref{theorem:dissipative-continuous}
shows that any energy profile $e(t)$ is indeed possible.

The existence of continuous solutions which violate the conservation of the total kinetic energy was first suggested in \cite{Onsager} by Onsager, where indeed he conjectured the existence of such solutions in $3$ space dimensions with any H\"older exponent  smaller than $\frac{1}{3}$. Onsager also asserted that such solutions do not exist if we impose the H\"older 
continuity with exponent larger than $\frac{1}{3}$ and this part of his conjecture was proved in \cite{Eyink} and
\cite{ConstantinETiti}. The considerations of Onsager are motivated by the Kolmogorov theory of isotropic
$3$-dimensional turbulence, where the phenomenon of anomalous dissipation in the Navier-Stokes equations
is postulated. This assumption seems to be widely confirmed experimentally, whereas no such
phenomenon is observed in $2$ dimensions. Indeed, for $d=2$ the conservation law for the enstrophy does
prevent it for solutions which start from sufficiently smooth initial data. However, the considerations put forward by Onsager which pertain the mathematical
structure of the equations do not depend on the dimension and this independence appears clearly also in
the proof of \cite{ConstantinETiti}, which works for any $d\geq 2$.

\medskip

The first proof of the existence of a weak solution violating the energy conservation was given in 
the groundbreaking work of Scheffer \cite{Scheffer93}, which showed the existence of a compactly supported nontrivial
weak solution in $\R^2\times \R$. A different construction of the existence of a compactly supported nontrivial
weak solution in $\mathbb T^2\times \R$ was then given by Shnirelman in \cite{Shnirelman1}. In both cases the 
solutions are only square summable as a function of both space and time variables. The first proof of the existence of a solution for which the total
kinetic energy is a monotone decreasing function has been
given by Shnirelman in \cite{Shnirelmandecrease}. Shnirelman's example is only in the energy space
$L^\infty ([0, \infty[, L^2 (\R^3))$. 

In the works \cite{DS1,DS2} these existence results were extended to solutions 
with bounded velocity and pressure and the same methods were also used to give quite 
severe counterexamples to the uniqueness
of admissible solutions, both for incompressible and compressible Euler. Further developments in fluid dynamics inspired by these works appeared subsequently in \cite{Chiodaroli, CFG, Shvydkoy, Szekelyhidi, SzWie, Wiedemann} and are surveyed in the note \cite{DSsurvey}. The paper \cite{DS1} 
introduced a new
point of view in the subject, highlighting connections to other counterintuitive solutions of 
(mainly geometric) systems of partial differential equations: in geometry these solutions are, according to Gromov, instances of the $h$-principle. 
It was also observed in \cite{CDSz} that the Onsager's Conjecture bears
striking similarities with the rigidity and flexibility properties of isometric embeddings of 
Riemannian manifolds, pioneered by the celebrated work of Nash \cite{Nash54}.
It turns out that the iteration procedure leading to Theorem \ref{theorem:dissipative-continuous} shares
some fundamental aspects with the approach of \cite{Nash54}: we refer to the introduction of \cite{DS3} for a thorough discussion.

\medskip

As already mentioned, we will also show that the 3D flows constructed in Theorem~\ref{theorem:dissipative-continuous} are genuinely three-dimensional. In order to formulate our statement precisely, consider a solution
$(v,p)$ of \eqref{eq:Euler} on $\mathbb T^3\times [0,T]$ and denote with the same letters
the corresponding solution on $\mathbb R^3\times [0,T]$ with the obvious periodicity in space.
The solution is then not genuinely $3$-dimensional if, after suitably changing coordinates in space, it takes the form
\begin{equation}\label{e:constv_3}
v (x,t) = (v_1 (x_1, x_2, t), v_2 (x_1, x_2, t), v_3)
\end{equation}
where $v_3$ is a constant. Now, a careful analysis of the proof of Theorem~\ref{theorem:dissipative-continuous}
will give the following statement.

\begin{theorem}\label{theorem:3D-flows-are-not-2D}
Let $d=3$. For any $\varepsilon>0$ the solutions constructed in Theorem~\ref{theorem:dissipative-continuous}
can be chosen to satisfy
\begin{equation}\label{e:large}
\sup_{t\in[0,1]}\left|\int_{\mathbb T^3}v(x,t)\otimes v(x,t)\,dx - \frac{1}{3}e(t) {\rm Id} \right|< \varepsilon,
\end{equation}
and
\begin{equation}
\sup_{t\in[0,1]}\|v\|_{H^{-1}(\mathbb T^3)}<\varepsilon.
\label{ineq:|v|H-1<epsilon}
\end{equation}
\end{theorem}

It is then easy to conclude that, for sufficiently small $\varepsilon$, a solution as in Theorem
\ref{theorem:3D-flows-are-not-2D} cannot have the form \eqref{e:constv_3}. Indeed, from
the bound \eqref{e:large} we would have $|v_3|\geq \frac{1}{6}\sqrt{\min_{t\in[0,1]}e(t)}$, whereas from
\eqref{ineq:|v|H-1<epsilon} we would have $|v_3|<\varepsilon$.

\noindent{\bf Remark}\qquad
A similar analysis in the two-dimensional case leads to the conclusion
that the two-dimensional flows constructed in Theorem~\ref{theorem:dissipative-continuous}
are typically not parallel flows.
However, it is classical that such flows are necessarily stationary,
and this is not possible if $e(t)$ is not constant.
\medskip

\section{Setup and general considerations on the construction}

As already mentioned, the proof of Theorem \ref{theorem:dissipative-continuous}
is based on an iteration procedure
in each step of which one solves a system of equations closely related to the Euler equations.
We introduce some terminology.
We let $\mathcal S^{d\times d}$ denote the space of symmetric $d\times d$ matrices
and $\mathcal S^{d\times d}_0$ its subspace of trace-free matrices.
Assume $v, p, \mathring R$ are smooth functions on $\mathbb T^d\times [0,T]$
taking values in $\mathbb R^d, \mathbb R$, and $\mathcal S^{d\times d}_0$ respectively.
We say that $(v,p,\mathring R)$ solves the Euler-Reynolds system if it satisfies
\begin{equation}
\partial_tv+{\rm div}\,(v\otimes v)+\nabla p={\rm div}\,\mathring R,\qquad {\rm div}\,v=0.
\label{eq:Euler-Reynolds}
\end{equation}

Theorem~\ref{theorem:dissipative-continuous} is a consequence of the following 
\begin{proposition}
\label{proposition:iteration}
Let $e\colon [0,T]\rightarrow \mathbb R$ be a smooth positive function.
Then, there exist constants $\eta$ and $M$, depending on $e$, with the following property.

Let $\delta\leq 1$ and $(v,p,\mathring R)$ be a solution of the Euler-Reynolds system
$$\partial_t v+{\rm div}\,(v\otimes v)+\nabla p={\rm div}\,\mathring R$$
such that
\begin{equation}
\frac34\delta e(t)\leq e(t)-\int_{\mathbb T^d} |v|^2(x,t)\,dx\leq \frac 54 \delta e(t),\qquad t\in[0,T]
\label{ineq:bound-on-e(t)-int|v|2}
\end{equation}
and 
\begin{equation}\sup_{x,t}|\mathring R(x,t)|\leq \eta\delta.
\label{ineq:bound-on-Ro}
\end{equation}
Then there exists $(v_1, p_1, \mathring R_1)$ 
solving the Euler-Reynolds system and satisfying the following estimates:
\begin{equation}
\frac38\delta e(t)\leq e(t)-\int |v_1|^2(x,t)\,dx\leq \frac58\delta e(t),\qquad t\in [0,T]
\label{ineq:bound-on-e(t)-int|v1|2}
\end{equation}
\begin{equation}
\sup_{x,t}|\mathring R_1(x,t)\leq \frac12\eta\delta,
\label{ineq:bound-on-R1o}
\end{equation}
\begin{equation}
\sup_{x,t}|v_1(x,t)-v(x,t)|\leq M\sqrt \delta,
\label{ineq:bound-on-v1-v}
\end{equation}
and
\begin{equation}
\sup_{x,t}|p_1(x,t)-p(x,t)|\leq M\delta.
\label{ineq:bound-on-p1-p}
\end{equation}
\end{proposition}

{\bf Proof of Theorem~\ref{theorem:dissipative-continuous}}:
\qquad 
Start with $v_0=0$, $p_0=0$, $\mathring R_0=0$, and $\delta=1$.
Apply Proposition~\ref{proposition:iteration} iteratively to obtain 
sequences $v_n$, $p_n$, and $\mathring R_n$ solving the Euler-Reynolds system (\ref{eq:Euler-Reynolds})
and satisfying
\begin{eqnarray}
\frac34\frac{e(t)}{2^n}\leq e(t)-\int_{\mathbb T^d}|v_n|^2(x,t)\,dx&\leq& \frac54\frac{e(t)}{2^n},\qquad t\in [0,T],\notag\\
\sup_{x,t}|\mathring R_n(x,t)|&\leq& \frac\eta{2^n},\notag\\
\sup_{x,t}|v_{n+1}(x,t)-v_n(x,t)|&\leq&M\sqrt{\frac1{2^n}},\notag\\
\sup_{x,t}|p_{n+1}(x,t)-p_n(x,t)|&\leq&\frac M{2^n}.\notag
\end{eqnarray}
The sequences $v_n$ and $p_n$ are Cauchy in $C(\mathbb T^d\times [0,T])$ and 
hence converge (uniformly) to continuous functions $v$ and $p$ respectively.
Likewise, $\mathring R_n$ converges (uniformly) to $0$.
Moreover, taking limits in the estimates on the energy, 
$$\int_{\mathbb T^d}|v|^2(x,t)\,dx=e(t),\qquad t\in[0,T].
$$
Also, we may pass to the limit in the (weak formulation of the) Euler-Reynolds system (\ref{eq:Euler-Reynolds})
and this shows that $v, p$ satisfy the (weak formulation of the) Euler equations (\ref{eq:Euler}).
\stopproof

In the rest of this Section we collect the main ingredients in the proofs of Theorems~\ref{theorem:dissipative-continuous}
and \ref{theorem:3D-flows-are-not-2D}.

\subsection{Construction of $(v_1,p_1,\mathring R_1)$}
Assume $(v,p,\mathring R)$ as in Proposition~\ref{proposition:iteration}.
Write $v_1=v+w$ and $p_1=p+q$ where $w$ and $q$ are to be determined.
Then,
\begin{eqnarray}
\partial_tv_1+{\rm div}\,(v_1\otimes v_1)+\nabla p_1
={\rm div}\,(w\otimes w+q{\rm Id}+\mathring R)
+\partial_tw+{\rm div}\,(v\otimes w+w\otimes v).\notag
\end{eqnarray}
The perturbation $w$ should be so chosen 
as to eliminate the first term in the right-hand side,
viewing the remainder as a small error.
Note how this first term is reminiscent of the stationary Euler equations.
Roughly speaking, the perturbation $w$ will be chosen as a high oscillation modulation
of a stationary solution.
More specifically, it will be taken of the form
$$w=w_o+w_c$$
where $w_o$ is a highly oscillator term given quite explicitly and 
the corrector term $w_c$ enforce the divergence-free condition ${\rm div}\,w=0$.
The oscillation term $w_o$ will depend on two parameters $\lambda$ and $\mu$ 
such that
$$\lambda, \mu, \frac\lambda\mu\in\mathbb N.$$
In fact, $\lambda$ and $\mu$ will have to be chosen sufficiently large,
depending on $v$, $\mathring R$, and $e$, 
in order that the desired estimates (\ref{ineq:bound-on-e(t)-int|v1|2}), (\ref{ineq:bound-on-R1o}),
(\ref{ineq:bound-on-v1-v}), and (\ref{ineq:bound-on-p1-p}) be satisfied.

\subsection{A linear set of stationary flows}
The stationary Euler equation is nonlinear.
The following Proposition says
that there exists a {\it linear} space of stationary solutions.
We first introduce some notation.
\begin{itemize}
\item {\bf Case $d=2$}\qquad
For $k\in \mathbb Z^2$, we let
\begin{equation}
b_k(\xi):=i\frac{k^\perp}{|k|}e^{ik\cdot\xi},
\qquad
\psi_k(\xi)=\frac{e^{ik\cdot\xi}}{|k|}%=-\frac{\cos(k\cdot \xi)}{|k|}.
\label{def:bk(xi)-in-2D}
\end{equation}
so that
$$b_k(\xi)=\nabla_\xi^\perp \psi_k(\xi),\quad{\rm where}\quad \nabla_\xi^\perp=(-\partial_{\xi^2},\partial_{\xi^1})
$$
hence ${\rm div}_\xi~b_k=0$.
Observe also that
$\overline{b_k}=b_{-k}$ and $\overline{\psi_k}=\psi_{-k}$.
\item {\bf Case $d=3$}\qquad 
For $k\in\mathbb Z^3$, we let
\begin{equation}
b_k(\xi):=B_ke^{ik\cdot\xi}
\label{def:bk(xi)-in-3D}
\end{equation}
where $B_k$ will be chosen appropriately, see Lemma~\ref{lemma:Beltrami-flows}.
In particular, one should have ${\rm div}_\xi~b_k=0$ and $\overline {b_k}=b_{-k}$.
\end{itemize}

\begin{lemma}[A linear set of stationary flows in 2D]
\label{lemma:steady-flows-in-2D}
Let $\nu\geq 1$.
For $k\in\mathbb Z^2$ such that $|k|^2=\nu$, let $a_k\in\mathbb C$ such that $\overline{a_k}=a_{-k}$.
Then,
$$
W(\xi)=\sum_{|k|^2=\nu}a_kb_k(\xi),
\qquad \Psi(\xi)=\sum_{|k|^2=\nu}a_k\psi_k(\xi)
$$
are $\mathbb R$-valued and satisfy
\begin{equation}
{\rm div}\,_\xi(W\otimes W)=\nabla_\xi\left(\frac{|W|^2}2-\nu\frac{\Psi^2}2\right).
\label{eq:div(WoW)=nabla(P)}
\end{equation}
Furthermore, 
\begin{equation}
\langle W\otimes W\rangle_\xi := \fint_{\mathbb T^2}W\otimes W\,d\xi=\sum_{|k|^2=\nu}|a_k|^2\left({\rm Id}-\frac k{|k|}\otimes \frac k{|k|}\right)
\label{eq:<WoW>}
\end{equation}
\end{lemma}
In other words, $W$ and $\Psi$ are the velocity field and stream function,
respectively, of a stationary flow with pressure $P=-\frac{|W|^2}2+\nu\frac{\Psi^2}2$.

\startproof
By direct computation one finds $\Delta_\xi\psi_k=-|k|^2\psi_k$, 
and hence that $\Delta_\xi \Psi=-\nu\Psi$.
Recall the identities
$${\rm div}\,_\xi(W\otimes W)=\frac12\nabla_\xi|W|^2+({\rm curl}_\xi W)W^\perp$$
where ${\rm curl}_\xi W=\partial_{\xi^1}W^2-\partial_{\xi^2}W^1=\Delta_\xi\Psi$ and $W^\perp=(-W^2,W^1)$.
Then,
$${\rm div}\,_\xi(W\otimes W)=\nabla_\xi\frac{|W|^2}2-\nu\Psi\nabla_\xi\Psi$$
as desired.

As for the average, write
\begin{eqnarray}
\fint_{\mathbb T^2}W\otimes  W(\xi)\,d\xi
&=&\sum_{j,k}a_ja_kb_j\otimes b_k\notag\\
&=&-\sum_{j,k}a_ja_ke^{i(k+j)\cdot \xi}\frac{j^\perp}{|j|}\otimes \frac{k^\perp}{|k|}\notag\\
&=&\sum_{j,k}a_k\overline{a_j}e^ {i(k-j)\cdot \xi}\frac{j^\perp}{|j|}\otimes \frac{k^\perp}{|k|}.
%&=&\sum_{j,k}a_ja_k \fint_{\mathbb T^2}\sin(j\cdot\xi)\sin(k\cdot \xi)\,d\xi \frac{j^\perp}{|j|}\otimes \frac{k^\perp}{|k|}.\notag
\end{eqnarray}
If $j\neq k$, $\fint_{\mathbb T^2}e^{i(k-j)\cdot\xi}\,d\xi=0$ and it is $1$ if $j=k$.
Thus,
$$
\langle W\otimes W \rangle_\xi=\sum_{|k|^2=\nu}|a_k|^2\frac{k^\perp}{|k|}\otimes \frac{k^\perp}{|k|}
=\sum_{|k|^2=\nu}|a_k|^2\left({\rm Id}-\frac k{|k|}\otimes \frac k{|k|}\right)
$$
where the last identity follows from $\frac{k^\perp}{|k|}\otimes \frac{k^\perp}{|k|}=\left({\rm Id}-\frac k{|k|}\otimes \frac k{|k|}\right)$
by direct calculation.
\stopproof

The following three-dimensional version was proved in \cite{DS3}.
\begin{lemma}[Beltrami flows]\label{lemma:Beltrami-flows}
Let $\nu\geq 1$.
There exist $B_k\in \mathbb C^3$ for $|k|=\nu$ such that, for any choice of $a_k\in\mathbb C$ such that $\overline{a_k}=a_{-k}$,
the vector field 
\begin{equation}
W(\xi)=\sum_{|k|=\nu}a_kb_k(\xi)
\label{def:W-Bk}
\end{equation}
is divergence-free and satisfies
\begin{equation}
{\rm div}\,\left(W\otimes W\right)=\nabla\frac{|W|^2}2.
\label{eq:div(WoW)=nabla|W|2/2}
\end{equation}
Furthermore,
\begin{equation}
\fint_{\mathbb T^3}W\otimes W\,d\xi=\sum_{|k|=\nu}|a_k|^2\left({\rm Id}-\frac k{|k|}\otimes \frac k{|k|}\right).
\label{eq:ave(WoW)=sum|ak|2Mk}
\end{equation}
\end{lemma}
Note that in the three-dimensional case the pressure is given by $P=-\frac{|W|^2}2$.

\startproof
This is Proposition~3.1 in \cite{DS3}.
\stopproof

\subsection{The geometric lemma}
\begin{lemma}[Geometric Lemma]
\label{lemma:geometric}
For every $N\in\mathbb N$ we can choose $r_0>0$ and $\nu\geq 1$ with the following property.
There exist pairwise disjoint subsets
$$
\Lambda_j\subset 
\left\{\begin{array}{ll}\{k\in \mathbb Z^2~:~|k|^2=\nu\}&{\rm if}\quad d=2\\\{k\in \mathbb Z^3~:~|k|^2=\nu\}&{\rm if}\quad d=3\end{array}\right.,
\qquad j\in\{1,\dots, N\}
$$
and smooth positive functions
$$
\gamma_k^{(j)}\in C^\infty(B_{r_0}({\rm Id})),\qquad j\in \{1, \dots, N\}, \quad k\in \Lambda_j
$$
such that
\begin{enumerate}
\item $k\in\Lambda_j$ implies $-k\in\Lambda_j$ and $\gamma_k^{(j)}=\gamma_{-k}^{(j)}$;
\item for each $R\in B_{r_0}({\rm Id})$ we have the identity
\begin{equation}
R=\sum_{k\in\Lambda_j}\left(\gamma_k^{(j)}(R)\right)^2\left({\rm Id}-\frac k{|k|}\otimes \frac k{|k|}\right)
\qquad R\in B_{r_0}({\rm Id}).
\label{eq:R=sum-of-Mk}
\end{equation}
\end{enumerate}
\end{lemma}
\startproof
This was proved for the case $d=3$ in Lemma~3.2 of \cite{DS3}
(up to a normalization constant).
Here we treat the case $d=2$ only.
For each $v\in \mathbb R^2\setminus\left\{0\right\}$ we define
$$M_v={\rm Id}-\frac v{|v|}\otimes\frac v{|v|}.$$

\noindent{\bf Step~1}\quad Fix $\nu\geq 1$ and for each set $F\subset\{k\in\mathbb Z^2~:~|k|^2=\nu\}$
denote $c(F)$ the interior of the convex hull in $\mathcal S^{2\times 2}$ (the space of symmetric $2\times 2$ matrices)
of the set $M_F:=\{M_k~:~k\in F\}$.
We claim in this step that it sufices to find $\nu$ and $N$ disjoint subsets $F_j\subset \{k\in\mathbb Z^2~:~ |k|^2=\nu\}$
such that
\begin{itemize}
\item $-F_j=F_j$;
\item $c(F_j)$ contains a positive multiple of the identity.
\end{itemize}
Indeed, we will show below that, if $F_j$ satisfies these two conditions,
then we can find $r_0>0$, a subset $\Gamma_j\subset F_j$, 
and smooth positive functions $\lambda_k^{(j)}\in C^\infty(B_{2r_0}({\rm Id}))$ for $k\in\Gamma_j$
such that
$$R=\sum_{k\in\Gamma_j}\lambda_k^{(j)}(R)M_k.$$
We then define the sets $\Lambda_j$ and the functions $\gamma_k^{(j)}$ for $k\in\Lambda_j$ according to
\begin{itemize}
\item $\Lambda_j:=\Gamma_j\cup -\Gamma_j$;
\item $\lambda_k^{(j)}=0$ if $k\in \Lambda_j\setminus \Gamma_j$;
\item $\gamma_k^{(j)}:=\sqrt{\lambda_k^{(j)}+\lambda_k^{(j)}}$ for $k\in\Lambda_j$.
\end{itemize}
Note that the sets $\Lambda_j$ are then symmetric and that the functions $\gamma_k^{(j)}$ satisfy (\ref{eq:R=sum-of-Mk}).
Moreover, since at least one of $\lambda_{\pm k}^{(j)}$ is positive on $B_{2r_0}({\rm Id})$,
$\gamma_k^{(j)}$ is smooth (and positive) in $B_{r_0}({\rm Id})$.

We now come to the existence of $\Gamma_j$.
The open set $c(F_j)$ contains an element $\alpha{\rm Id}$ for some $\alpha>0$.
Observe that the space $\mathcal S^{2\times 2}$ of symmetric $2\times 2$ matrices 
has dimension $3$.
Since $\alpha{\rm Id}$ sits inside the open set $c(F_j)$, there exists a $4$-simplex $S$
with vertices $A_1, \dots, A_4\in c(F_j)$
such that $\alpha{\rm Id}$ belongs to the interior of $S$.
Let then $\vartheta$ so that the ball $\tilde U$ centered at $\alpha{\rm Id}$ and radius $\vartheta$ si contained in $S$.
Then, each $R\in\tilde U$ can be written in a unique way as a convex combination of $A_i$'s:
$$R=\sum_{i=1}^4\beta_i(R)A_i$$
where the functions $\beta_i$ are positive and smooth in $\tilde U$.

Using now Carath\'eodory's Theorem, each $A_i$ is the convex combination $\sum_n \lambda_{i,n}M_{v_{i,n}}$
of at most $4$ matrices $M_{v_{i,n}}$ where $v_{i,n}\in F_j$,
where we require that each $\lambda_{i,n}$ be positive.
(Carath\'eodory's Theorem guarantees the existence of $4$ elements $M_{v_{i,n}}$ such that $A_i$ 
belongs to their {\it closed} convex hull.
If we insist that all coefficients should be positive, 
then some of these elements should be thrown away.)

Set now $r_0:=\frac\vartheta{2\alpha}$.
Then
$$R=\sum_{i,n}\frac1\alpha\beta_i(\alpha R)\lambda_{i,n}M_{v_{i,n}},\qquad R\in B_{2r_0}({\rm Id})$$
and each coefficient $\frac1\alpha\beta_i(\alpha R)\lambda_{i,n}$ is positive for $R\in B_{2r_0}({\rm Id})$.
The set $\Gamma_j$ is then taken as $\{v_{i,n}\}$.
Since one might have $v_{i,n}=v_{l,m}$ for distinct $(i,n)$ and $(l,m)$, the function $\lambda_k$ will be defined by
$$\lambda_k(R)=\sum_{(i,n):v_{i,n}=k}\frac1\alpha\beta_i(\alpha R)\lambda_{i,n}$$
and this completes Step~1.
\bigskip

\noindent{\bf Step~2}\quad By Step~1, in order to prove the lemma,
it suffices to find $\nu$ and $N$ disjoint families $F_1, \cdots, F_N\subset \sqrt\nu\mathbb S^1\cap \mathbb Z^2$
such that each set $c(F_j)$ contains a positive multiple of the identity.
Note that $\mathbb S^1\cap \mathbb Q^2$ is dense in $\mathbb S^1$.
Indeed, let
$$u\in\mathbb R\mapsto s(u):=\left(\frac{2u}{u^2+1},\frac{u^2-1}{u^2+1}\right)\in\mathbb S^1\subset\mathbb R^2.$$
Clearly $s(\mathbb Q)\subset \mathbb Q^2$. Since $\mathbb Q$ is dense in $\mathbb R$
and $s$ is a diffeomorphism onto $\mathbb S^1\setminus\{(0,1)\}$, the claim is proved.

In turn, there exists a sequence $\nu_k\rightarrow +\infty$ such that the sets $\mathbb S^1\cap \frac1{\sqrt{\nu_k}}\mathbb Z^2$
converge, in the Hausdorff sense, to the entire circle $\mathbb S^1$.
Given the sequence $\nu_k$, one can easily partition each $\sqrt{\nu_k}\mathbb S^1\cap\mathbb Z^2$ into $N$ 
disjoint symmetric families $\{F_j^k\}_{j=1,\dots,N}$ in such a way that,
for each fixed $j$, the corresponding sequence of sets $\left\{\frac1{\sqrt{\nu_k}}F_j^k\right\}_k$
converges in the Hausdorff sense to $\mathbb S^1$.
Hence, any point of $c(\mathbb S^1)$ is contained in $c(\frac1{\sqrt{\nu_k}}F_j^k)$ for $k$ sufficiently large.
On the other hand,
it is easy to see that $c(\mathbb S^1)$ contains a multiple $\alpha{\rm Id}$ of the identity.
(One can adapt for instance the argument of Lemma~4.2 in \cite{DS1}.)

By Step~1, this concludes the proof.
\stopproof

\section{The maps $v_1, p_1,\mathring R_1$}
Let $e(t)$ be as in Theorem~\ref{theorem:dissipative-continuous}
and suppose $(v,p,\mathring R)$ satisfies the hypothesis of Proposition~\ref{proposition:iteration}.
In this Section we define the next iterates $v_1=v+w$, $p_1=p+q$, and $\mathring R_1$.
The perturbation $w$ will be defined as $w=w_o+w_c$ where $w_o$ is given by an explicit formula,
see (\ref{eq:wo(x,t)}),
and the corrector $w_c$ guarantees that ${\rm div}\,w=0$.

\subsection{The perturbation $w_o$}
We apply Lemma~\ref{lemma:geometric} with 
$$N=2^d$$
(where $d\in\{2,3\}$ denotes the number of space dimensions)
to obtain $\nu\geq 1$ and $r_0>0$,
and pairwise disjoint families $\Lambda_j$ with corresponding functions $\gamma_k^{(j)}\in C^\infty(B_{r_0}({\rm Id}))$.

The following is proved in \cite{DS3}
and works for any number of space dimensions.
Denote $\mathcal C_1, \dots, \mathcal C_N$ the equivalence classes 
of $\mathbb Z^d/\sim$ where $k\sim l$ if $k-l\in (2\mathbb Z)^d$.
The parameter $\mu\in \mathbb N$ will be fixed later.
\begin{proposition}[Partition of the space of velocities]
There exists a partition of the space of velocities, 
namely $\mathbb R$-valued functions $\alpha_l(v)$ for $l\in\mathbb Z^d$ satisfying
\begin{equation}
\sum_{l\in\mathbb Z^d}(\alpha_l(v))^2\equiv 1
\label{eq:sum-of-ak-squared=1}
\end{equation}
such that, setting
\begin{equation}
\phi_k^{(j)}(v,\tau)=\sum_{l\in\mathcal C_j}\alpha_l(\mu v)e^{-i(k\cdot\frac l\mu)\tau},\qquad j=1,\dots, N,\quad k\in\mathbb Z^d,
\label{eq:phi(j)k}
\end{equation}
we have the following estimates:
\begin{eqnarray}
\sup_{v,\tau}|D^m_v\phi_k^{(j)}(v,\tau)|&\leq& C(m,d)\mu^m,\label{ineq:estimates-v-derivatives-phi(j)k}\\
\sup_{v,\tau}|D^m_v(\partial_\tau\phi_k^{(j)}+i(k\cdot v)\phi_k^{(j)})(v,\tau)|&\leq& C(m,|k|,d)\mu^{m-1}.\label{ineq:estimates-transport-phi(j)k}
\end{eqnarray}
Furthermore, $\overline{\phi_k^{(j)}}=\phi_{-k}^{(j)}$ and
\begin{equation}|\phi_k^{(j)}(v,\tau)|^2=\sum_{l\in\mathcal C_j}\alpha_l^2(v).
\label{eq:|phi(j)k|2=sum-of-alpha(k)2}
\end{equation}
\end{proposition}

Set now
\begin{equation}
\rho(s):=\frac1{d(2\pi)^d}\left(e(t)\left(1-\frac\delta2\right)-\int_{\mathbb T^d}|v|^2(x,t)\,dx\right)
\label{eq:rho(s)}
\end{equation}
and
\begin{equation}
R(y,s):=\rho(s){\rm Id}-\mathring R(y,s)
\label{eq:R(y,s)}
\end{equation}
(observe that the tensor $-\mathring{R}$ is then the traceless part of $R$).
Define
\begin{equation}
w_o(x,t):=W(x,t,\lambda t, \lambda x)
%=\sqrt{\rho(t)}\sum_{j=1}^4\sum_{k\in\Lambda_j}\gamma_k^{(j)}\left(\frac{R(x,t)}{\rho(t)}\right)\phi_k^{(j)}(v(x,t),\lambda t)b_k(\lambda x)
\label{eq:wo(x,t)}
\end{equation}
where
\begin{eqnarray}
W(y,s,\tau,\xi)
&=&\sum_{|k|^2=\nu}a_k(y,s,\tau)b_k(\xi)\notag\\
&=&\sqrt{\rho(s)}\sum_{j=1}^N\sum_{k\in\Lambda_j}\gamma_k^{(j)}\left(\frac{R(y,s)}{\rho(s)}\right)\phi_k^{(j)}(v(y,s),\tau)b_k(\xi).
\label{eq:W(y,s,tau,xi)}
\end{eqnarray}
(The velocity fields $b_k$ were defined in (\ref{def:bk(xi)-in-2D}) for $d=2$ and in
\eqref{def:bk(xi)-in-3D} for $d=3$.)
For $d=2$ we also introduce the corresponding stream function
\begin{equation}
\psi_o(x,t):=\Psi(x,t,\lambda t, \lambda x)
%=\sqrt{\rho(t)}\sum_{j=1}^4\sum_{k\in\Lambda_j}\gamma_k^{(j)}\left(\frac{R(x,t)}{\rho(t)}\right)\phi_k^{(j)}(v(x,t),\lambda t)\psi_k(\lambda x)
\label{eq:psio(x,t)}
\end{equation}
where
\begin{eqnarray}
\Psi(y,s,\tau,\xi)
&=&\sum_{|k|^2=\nu}a_k(y,s,\tau)\psi_k(\xi)\notag\\
&=&\sqrt{\rho(s)}\sum_{j=1}^N\sum_{k\in\Lambda_j}\gamma_k^{(j)}\left(\frac{R(y,s)}{\rho(s)}\right)\phi_k^{(j)}(v(y,s),\tau)\psi_k(\xi).\notag
\end{eqnarray}
(The stream functions $\psi_k$ were defined in (\ref{def:bk(xi)-in-2D}).)

\subsection{The constants $\eta$ and $M$}
In this Section we fix the values of the constants $\eta$ and $M$ from Proposition~\ref{proposition:iteration}.

The perturbation $w_o$ is well defined provided $\frac R\rho\in B_{r_0}({\rm Id})$
where $r_0$ is as in Lemma~\ref{lemma:geometric}.
By definition (\ref{eq:rho(s)}) of $\rho$ and assumption (\ref{ineq:bound-on-e(t)-int|v|2}),
$$\rho(t)\geq \frac1{d(2\pi)^d}\frac\delta4e(t)\geq c\delta \min_{t\in [0,T]}e(t)$$
where $c=\frac1{4d(2\pi)^d}$.
Then,
$$\left|\frac R{\rho(t)}-{\rm Id}\right|\leq \frac1{c\delta \min_{t\in [0,T]}e(t)}\left|\mathring R\right|\leq \frac\eta{c \min_{t\in [0,T]}e(t)}.$$
Thus, we choose $\eta$ satisfying
\begin{equation}\eta\leq \frac12c\min_{t\in[0,T]}e(t)r_0.\label{def:eta}
\end{equation}
Observe that this restriction is independent of $\delta$.

We choose first a constant $M'>1$ such that
\begin{equation}\label{def:M'}
M'>
4\left(\sum_{j=1}^N\sum_{k\in\Lambda_j}\,\sup_R|\gamma_k^{(j)}(R)|\,\sup_{v,\tau}|\phi_k^{(j)}(v,\tau)|\,\sup_\xi|b_k(\xi)|\right)^2
\end{equation}
%and, in case $d=2$, we additionally impose that
%\begin{equation}
%M'>\frac1{\sqrt\nu} \sum_{j=1}^N\sum_{k\in\Lambda_j}\,\sup_R|\gamma_k^{(j)}(R)|\,\sup_{v,\tau}|\phi_k^{(j)}(v,\tau)|\,\sup_\xi|\psi_k(\xi)|
%\label{def:M'}
%\end{equation}
and then choose $M>1$ such that 
$$M\geq M' \frac{\sup_{0\leq t\leq T}e(t)}{4d(2\pi)^d}.$$
Observe that with these choices we have
\begin{equation}
\max\big\{\sqrt\nu\|\psi_o\|_0, \|w_o\|_0\big\}\quad\leq\quad \frac{\sqrt M}2\sqrt\delta.
\label{ineq:fix-parameter-M}
\end{equation}
since $\rho(t)\leq \frac1{d(2\pi)^d}\delta e(t)$ by (\ref{eq:rho(s)}) and (\ref{ineq:bound-on-e(t)-int|v|2}).

\subsection{The correction $w_c$}
To obtain $w$ from $w_o$ we need to introduce the Leray projection onto divergence-free vector fields with zero average.
\begin{definition}[The Leray projector]
Let $v\in C^\infty(\mathbb T^d,\mathbb R^d)$ be a smooth vector field.
Let
$$
\mathcal Qv:=\nabla \phi+\fint_{\mathbb T^d}v
$$
where $\phi\in C^\infty(\mathbb T^d)$ is the solution to
$$
\Delta\phi={\rm div}\,v\quad {\rm in}~\mathbb T^d\qquad{\rm subject~to}\quad \fint_{\mathbb T^d}\phi=0.
$$
We denote by
$$
\mathcal P:=I-\mathcal Q
$$
the Leray projector onto divergence-free vector fields with zero average.
\end{definition}

The iterate $v_1$ is then expressed as
\begin{equation}
v_1=v+\mathcal Pw_o=v+w_o+w_c,\qquad w_c=-\mathcal Qw_o=w-w_o.
\label{eq:v1=v+wo+wc}
\end{equation}

\subsection{The pressure term $p_1$}
We set
\begin{equation}
p_1:=\left\{\begin{array}{ll}p-\left(\frac{|w_o|^2}2-\nu\frac{\psi_o^2}2\right)&{\rm if}\quad d=2,\\
p-\frac{|w_o|^2}2&{\rm if}\quad d=3.\end{array}\right.
\label{eq:p1}
\end{equation}
This choice for $p_1$ will become clearer at the end of the proof, see Lemma~\ref{lemma:estimate-on-oscillation-part}.

\subsection{The tensor $\mathring R_1$}
To construct $\mathring R_1$, we introduce another operator.
\begin{definition}
Let $v\in C^\infty(\mathbb T^d,\mathbb R^d)$ be a smooth vector field.
We define $\mathcal Rv$ to be the matrix-valued periodic function
\begin{equation}
\mathcal Rv=\left\{\begin{array}{ll}
\nabla u+(\nabla u)^\top -({\rm div}\,\,u){\rm Id}&{\rm if}\quad d=2\\
\frac14\left(\nabla\mathcal Pu+(\nabla\mathcal Pu)^\perp\right)+\frac34\left(\nabla u+(\nabla u)^\perp\right)-\frac12({\rm div}\,u){\rm Id}&{\rm if}\quad d=3
\end{array}\right.
\label{eq:R=divInverse}
\end{equation}
where $u\in C^\infty(\mathbb T^d,\mathbb R^d)$ is the solution to
$$
\Delta u=v-\fint_{\mathbb T^d}v\quad {\rm in}~\mathbb T^d,\qquad {\rm subject~to}\quad \fint_{\mathbb T^d}u=0.
$$
\end{definition}
The operator $\mathcal R$ satisfies the following properties.
\begin{lemma}[$\mathcal R={\rm div}\,^{-1}$]\label{lemma:R=divInverse}
For any $v\in C^\infty(\mathbb T^d,\mathbb R^d)$ we have
\begin{enumerate}
\item $\mathcal Rv(x)$ is a symmetric trace-free matrix for each $x\in\mathbb T^d$;
\item ${\rm div}\,\mathcal Rv=v-\fint_{\mathbb T^d}v$.
\end{enumerate}
\end{lemma}
\startproof
The case $d=3$ is treated in \cite{DS3} and thus we assume $d=2$.
Clearly, $\mathcal Rv$ is symmetric.
Next,
$$
{\rm tr}\,\mathcal Rv=2{\rm div}\,u-2{\rm div}\,u=0
$$
and
\begin{eqnarray}
{\rm div}\,\mathcal Rv
&=&\Delta u+{\rm div}\,(\nabla u)^\top-{\rm div}\,\left(({\rm div}\,u){\rm Id}\right)\notag\\
&=&\Delta u+{\rm div}\,(\nabla u)^\top-\nabla {\rm div}\,u\notag\\
&=&\Delta u\notag\\
&=&v-\fint_{\mathbb T^2}v\notag
\end{eqnarray}
by definition of $u$.
\stopproof

Then we set
\begin{equation}
\mathring R_1:=\mathcal R\left(\partial_tv_1+{\rm div}\,(v_1\otimes v_1)+\nabla p_1\right).
\label{eq:R1o}
\end{equation}
One verifies, as in \cite{DS3} (after Lemma~4.3, p.~15),
that the argument in the right-hand side 
has zero average:
${\rm div}\,(v_1\otimes v_1+p_1{\rm Id})$ clearly has average zero,
and so does $\partial_tv_1=\partial_tv+\partial_tw$ 
since $\partial_tv=-{\rm div}\,(v\otimes v+p{\rm Id})$ has average zero as well as $w$ by definition of $\mathcal P$.
In turn, Lemma~\ref{lemma:R=divInverse} yields
$$
\partial_tv_1+{\rm div}\,(v_1\otimes v_1)+\nabla p_1={\rm div}\,\mathring R_1.
$$

\section{Estimates}

The letter $m$ will denote a natural number (in $\mathbb N$), and $\alpha$ a real number in the interval $(0,1)$.
The letter $C$ will always denote a generic constant which may depend on $e,v,\mathring R, \nu, \alpha$, and $\delta$,
but not on $\lambda$ nor $\mu$.
We will further impose that 
\begin{equation}1\leq \mu\leq \lambda.
\label{ineq:1<mu<lambda}
\end{equation}

The sup-norm is denoted $\|f\|_0=\sup_{\mathbb T^d}|f|$.
The H\"older seminorms are given by
\begin{eqnarray}
\left[f\right]_m&:=&\max_{|\gamma|=m}\|D^\gamma f\|_0,\notag\\
\left[f\right]_{m+\alpha}&:=&\max_{|\gamma|=m}\sup_{x\neq y}\frac{|D^\gamma f(x)-D^\gamma f(y)|}{|x-y|^\alpha}\notag
\end{eqnarray}
and the H\"older norms are given by
\begin{eqnarray}
\|f\|_m&:=&\sum_{j=0}^m\left[f\right]_j,\notag\\
\|f\|_{m+\alpha}&:=&\|f\|_m+\left[f\right]_{m+\alpha}.\notag
\end{eqnarray}
We also recall the following elementariy identity:
for $0\leq r\leq 1$
$$\left[fg\right]_r\leq C(\left[f\right]_r\|g\|_0+\|f\|_0\left[g\right]_r).$$

\subsection{Schauder estimates}

The next Proposition collects estimates in H\"older spaces (``Schauder estimates'')
for various operators used in the remainder.
\begin{proposition}
For any $\alpha\in (0,1)$ and $m\in\mathbb N$ there exists a constant $C=C(m,\alpha,d)$ satisfying the following properties.
If $\phi,\psi\colon\mathbb T^d\rightarrow \mathbb R$ are the unique solutions to
$$
\left\{\begin{array}c\Delta\phi=f\\\fint_{\mathbb T^d} \phi=0\end{array}\right.,\qquad 
\left\{\begin{array}c\Delta\phi={\rm div}\,F\\\fint_{\mathbb T^d} \psi=0\end{array}\right.,\qquad 
$$
then
$$\|\phi\|_{m+2,\alpha}\leq C\|f\|_{m,\alpha},\quad {\rm and}\quad \|\psi\|_{m+1,\alpha}\leq C\|F\|_{m,\alpha}.$$
Moreover, we have the following estimates:
\begin{eqnarray}
\|\mathcal Qv\|_{m+\alpha}&\leq&C\|v\|_{m+\alpha},\label{ineq:Schauder-estimate-for-Q}\\
\|\mathcal Pv\|_{m+\alpha}&\leq&C\|v\|_{m+\alpha},\label{ineq:Schauder-estimate-for-P}\\
\|\mathcal Rv\|_{m+1,\alpha}&\leq&C\|v\|_{m+\alpha},\label{ineq:Schauder-estimate-for-R}\\
\|\mathcal R{\rm div}\,A\|_{m+\alpha}&\leq&C\|A\|_{m+\alpha},\label{ineq:Schauder-estimate-for-Rdiv}\\
\|\mathcal R\mathcal Q{\rm div}\,A\|_{m+\alpha}&\leq&C\|A\|_{m+\alpha}.\label{ineq:Schauder-estimate-for-RQdiv}
\end{eqnarray}
\end{proposition}
These estimates are fairly standard and a detailed proof is given in \cite{DS3} in three dimensions.
It is easy easy to see that they also hold in two dimensions as well: the only operator which is
defined differently is $\mathcal R$ and this makes the $2$-dimensional case only easier.
Suffice it to say that these estimates are the expected ones:
$\mathcal P$ and $\mathcal Q$ are differential operators of degree $0$;
$\mathcal R$ is a differential operator of degree $-1$;
and div is a differential operator of degree $1$.

The effect of the oscillation parameter $\lambda$ is described in the following
\begin{proposition}\label{proposition:oscillations-estimates}
Let $k\in\mathbb Z^d\setminus 0$ and $\lambda\geq 1$ be fixed.
\begin{enumerate}
\item For any $a\in C^\infty(\mathbb T^d)$ and $m\in\mathbb N$ we have
$$\left|\int_{\mathbb T^d}a(x)e^{i\lambda k\cdot x}\,dx\right|\leq \frac{\left[a\right]_m}{\lambda^m}.$$
\item Let $\phi_\lambda\in C^\infty(\mathbb T^d)$ be the solution to
$$\Delta\phi_\lambda=f_\lambda\quad {\rm in}\quad \mathbb T^d$$
subject to $\fint_{\mathbb T^d}\phi_\lambda=0$
where $f_\lambda(x)=a(x)e^{i\lambda k\cdot x}-\fint_{\mathbb T^d}a(y)e^{i\lambda k\cdot y}dy$.
Then, for any $\alpha \in (0,1)$ we have the estimate
$$
\|\nabla\phi_\lambda\|_\alpha 
\leq \frac C{\lambda^{1-\alpha}}\|a\|_0+\frac C{\lambda^{m-\alpha}}[a]_m+\frac C{\lambda^m}[a]_{m+\alpha}
$$
where $C=C(m,\alpha,d)$.
\end{enumerate}
\end{proposition}
This was established in \cite{DS3} and is in fact valid in any dimension.

The following is a consequence of the definition (\ref{eq:R=divInverse}) of $\mathcal R$, 
the Schauder estimate (\ref{ineq:Schauder-estimate-for-P}) for $\mathcal P$,
and Proposition~\ref{proposition:oscillations-estimates}.
\begin{corollary}[Estimates for the operator $\mathcal R$]\label{corollary:estimates-for-R}
Let $k\in\mathbb Z^d\setminus 0$ be fixed.
For a smooth vector field $a\in C^\infty(\mathbb T^d,\mathbb R^d)$, let $F(x):=a(x)e^{i\lambda k\cdot x}$.
Then, we have
$$
\|\mathcal R(F)\|_\alpha\leq \frac C{\lambda^{1-\alpha}}\|a\|_0+\frac C{\lambda^{m-\alpha}}[a]_m+\frac C{\lambda^m}[a]_{m+\alpha}
$$
for $C=C(m,\alpha,d)$.
\end{corollary}

\subsection{Estimates on the corrector and the energy}
We recall that $w_o(x,t)$ is defined in (\ref{eq:wo(x,t)}).
In what follows, if $f$ is a function of time and space, we will use
the notation
$\|f\|_r$ for the ``H\"older norm in space'', that is
\[
\|f\|_r = \sup_{t\in [0,T]} \|f (\cdot, t)\|_r\, .
\]

\begin{lemma}
\label{lemma:estimates-on-coefficients}
\begin{enumerate}
\item Let $a_k\in C^\infty(\mathbb T^d\times [0,T]\times \mathbb R)$ be as in (\ref{eq:wo(x,t)}).
Then, for any $r\geq 0$,
\begin{eqnarray}
\|a_k(\cdot, s,\tau)\|_r&\leq&C\mu^r,\notag\\
\|\partial_sa_k(\cdot,s,\tau)\|_r&\leq& C\mu^{r+1},\notag\\
\|\partial_\tau a_k(\cdot,s,\tau)\|_r&\leq& C\mu^r,\notag\\
\|(\partial_\tau a_k+i(k\cdot v)a_k)(\cdot,s,\tau)\|_r&\leq&C\mu^{r-1}.\notag
\end{eqnarray}
\item The matrix-valued function $W\otimes W$ is given by
\begin{eqnarray}
(W\otimes W)(y,s,\tau,\xi)
=R(y,s)
+\sum_{1\leq |k|\leq 2\nu}U_k(y,s,\tau)e^{ik\cdot\xi}
\label{eq:WoW=R+SumUk}
\end{eqnarray}
where the coefficients $U_k\in C^\infty(\mathbb T^d\times [0,T]\times \mathcal S^{d\times d})$ satisfy, for any $r\geq 0$,
\begin{eqnarray}
\|U_k(\cdot, s,\tau)\|_r&\leq&C\mu^r,\notag\\
\|\partial_sU_k(\cdot,s,\tau)\|_r&\leq& C\mu^{r+1},\notag\\
\|\partial_\tau U_k(\cdot,s,\tau)\|_r&\leq& C\mu^r,\notag\\
\|(\partial_\tau+i(k\cdot v)U_k)(\cdot,s,\tau)\|_r&\leq&C\mu^{r-1}.\notag
\end{eqnarray}
\item For $d=2$, we have
$$
\frac{|W|^2}2-\nu\frac{\Psi^2}2
=\sum_{1\leq |k|\leq 2\nu}\tilde a_ke^{ik\cdot \xi}
$$
where the functions $\tilde a_k$ satisfy for any $r\geq 0$
$$
\|\tilde a_k(\cdot, s,\tau)\|_r\leq C\mu^r.
$$
\end{enumerate}
In all these estimates the constant $C$ depends on $r$, $e$, $v$, and $\mathring R$,
but is independent of $(s,\tau)$ and $\mu$.
\end{lemma}
\startproof
The first two items are proved in \cite{DS3} for $d=3$
and we only briefly recall their proof since it is identical for $d=2$.
The estimates on the coefficients $a_k$ follow from the estimates (\ref{ineq:estimates-v-derivatives-phi(j)k}) and
(\ref{ineq:estimates-transport-phi(j)k}) on $\phi_k^{(j)}$.
Next, write $W\otimes W$ as a Fourier series in $\xi$:
$$W\otimes W(y,s,\tau,\xi)=U_0(y,s,\tau)+\sum_{1\leq |k|\leq 2\nu}U_k(y,s,\tau)e^{ik\cdot \xi}$$
where the entries in the $U_k$'s are quadratic in the $a_k$'s and $\|a_k\|_0\leq C$.
Thus, the $U_k$'s satisfy the claimed estimates.
Furthermore, 
\begin{eqnarray}
U_0(y,s,\tau)&=&\fint_{\mathbb T^d}W\otimes W\,d\xi\notag\\
&\stackrel{(\ref{eq:<WoW>})}=&\rho\sum_{j=1}^N\sum_{k\in\Lambda_j}\left(\gamma_k^{(j)}\left(\frac R\rho\right)\right)^2|\phi_k^{(j)}(v,\tau)|^2\left({\rm Id}-\frac k{|k|}\otimes \frac k{|k|}\right)\notag\\
&\stackrel{(\ref{eq:|phi(j)k|2=sum-of-alpha(k)2})}=&\rho\sum_{j=1}^N\sum_{k\in\Lambda_j}\sum_{l\in\mathcal C_j}\left(\gamma_k^{(j)}\left(\frac R\rho\right)\right)^2\alpha_l^2(v)\left({\rm Id}-\frac k{|k|}\otimes \frac k{|k|}\right)\notag\\
&\stackrel{(\ref{eq:R=sum-of-Mk})}=&R\sum_{j=1}^N\sum_{l\in\mathcal C_j}\alpha_l^2(v)\notag\\
&\stackrel{(\ref{eq:sum-of-ak-squared=1})}=&R.
\end{eqnarray}

As for the third item of the Proposition, concerning $\frac{|W|^2}2-\nu\frac{\Psi^2}2$ when $d=2$,
we compute, omitting variables and remembering that the sums are over $j$ and $k$ such that $|j|^2=|k|^2=\nu$,
\begin{eqnarray}
\Psi^2
&=&\sum_{j,k}a_ja_k\psi_j\psi_k\notag\\
&=&\frac1\nu\sum_{j,k}a_ja_ke^{i(k+j)\cdot\xi}\notag\\
&=&\frac1\nu\sum_{j,k}a_k\overline{a_j}e^{i(k-j)\cdot\xi}\notag\\
&=&\frac1\nu\sum_{|k|^2=\nu}|a_k|^2+\frac1\nu\sum_{j\neq k}a_k\overline{a_j}e^{i(k-j)\cdot\xi}\notag\\
&=&\frac1\nu\sum_{|k|^2=\nu}|a_k|^2+\sum_{1\leq |k|\leq 2\nu}\underline a_k(y)e^{i(k\cdot \xi)}\notag
\end{eqnarray}
where the coefficients $\underline a_k$ are quadratic in the $a_k$'s.
But from the expression (\ref{eq:<WoW>}) for $\langle W\otimes W\rangle_\xi$ we deduce
\begin{eqnarray}
|W|^2={\rm tr}\,(W\otimes W)
&=&{\rm tr}\,R+\sum_{1\leq |k|\leq 2\nu}{\rm tr}\,U_ke^{ik\cdot \xi}\notag\\
&=&\sum_{|k|^2=\nu}|a_k|^2+\sum_{1\leq |k|\leq 2\nu}{\rm tr}\,U_ke^{ik\cdot \xi}.\notag
\end{eqnarray}
Subtracting the above expression for $\nu\Psi^2$ from that of $|W|^2$ we obtain the desired expression
with suitable coefficients $\tilde a_k$ which satisfy the same estimates as the $a_k$'s.
\stopproof

\begin{lemma}[Estimate on the corrector]\label{lemma:estimate-on-the-corrector}
\begin{equation}
\|w_c\|_\alpha\leq C\frac\mu{\lambda^{1-\alpha}}.
\label{ineq:estimate-on-wc}
\end{equation}
\end{lemma}
\startproof
This is proved in Lemma~6.2, p.~21, in \cite{DS3} for $d=3$.
For clarity we reprove it for $d=2$ with the appropriate adjustments.
Recall that
$$
w_o(x,t)=\sum_{|k|^2=\nu} a_k(x,t,\lambda t)\nabla^\perp\psi_k(\lambda x)=\nabla_\xi^\perp \Psi(x,t,\lambda t, \lambda x)
$$
where $\Psi(y,s,\tau,\xi) =\sum_{|k|^2=\nu} a_k(y,s,\tau)\psi_k(\xi)$.
But
$$
\nabla^\perp\left(a_k(x,t,\lambda t)\psi_k(\lambda x)\right)
=\lambda a_k(x,t,\lambda t)(\nabla^\perp\psi_k)(\lambda x)+\psi_k(\lambda x)\nabla^\perp a_k(x,t,\lambda t)
$$
or
$$
a_k(x,t,\lambda t)(\nabla^\perp\psi_k)(\lambda x)
=\frac1\lambda\left\{\nabla^\perp\left(a_k(x,t,\lambda t)\psi_k(\lambda x)\right)-\psi_k(\lambda x)\nabla^\perp a_k(x,t,\lambda)\right\}
$$
and thus
\begin{eqnarray}
w_o(x,t)
%&=&\sum_k a_k(x,t,\lambda t)b_k(\lambda x)\notag\\
&=&\frac1\lambda \nabla^\perp\left(\sum_{|k|^2=\nu} a_k(x,t,\lambda t)\psi_k(\lambda x)\right)
-\frac1\lambda \sum_{|k|^2=\nu} \psi_k(\lambda x)\nabla^\perp a_k(x,t,\lambda t).\notag
\end{eqnarray}
Since $\mathcal Q$ eliminates the divergence-free part and ${\rm div}\,\circ \nabla^\perp=0$, 
we have by (\ref{eq:v1=v+wo+wc})
\begin{equation}
w_c(x,t)=\mathcal Qw_o(x,t)
=\frac1\lambda \mathcal Q\left(\sum_{|k|^2=\nu} \psi_k(\lambda x)\nabla^\perp a_k(x,t,\lambda t)\right)
=:\frac1\lambda \mathcal Qu_c(x,t).
\label{eq:introduce-uc}
\end{equation}
Thus, by Schauder estimate (\ref{ineq:Schauder-estimate-for-Q}) for $\mathcal Q$ 
and the estimates on the coefficients $a_k$ from Lemma~\ref{lemma:estimates-on-coefficients}, we find
$$
\|w_c\|_\alpha\leq \frac C\lambda \|u_c\|_\alpha\leq C\frac\mu{\lambda^{1-\alpha}}.
$$
\stopproof

\begin{lemma}[Estimate on the energy]\label{lemma:estimate-on-the-energy}
\begin{equation}
\left|e(t)(1-\frac\delta2)-\int_{\mathbb T^d}|v_1|^2\,dx\right|\leq C\frac\mu{\lambda^{1-\alpha}}.
\label{ineq:estimate-on-e(t)-int|v1|2}
\end{equation}
\end{lemma}
\startproof
This is proved in Lemma~6.3, p.~21 of \cite{DS3} for $d=3$.
For clarity we briefly recall it in the case $d=2$.
Taking the trace in the expression (\ref{eq:WoW=R+SumUk}) for $W\otimes W$ gives
$$|W|^2={\rm tr}\,R+\sum_{1\leq |k|\leq 2\nu}{\rm tr}\,U_ke^{ik\cdot \xi}.$$
From part (1) of Proposition~\ref{proposition:oscillations-estimates} with $m=1$, 
and using estimates on $U_k$ from Lemma~\ref{lemma:estimates-on-coefficients}, we find
$$\left|\int_{\mathbb T^d}|w_o|^2-{\rm tr}\,\,R\,dx\right|\quad{\rm and}\quad
\left|\int_{\mathbb T^d}v\cdot w_o\,dx\right|\quad \leq\quad C\frac\mu\lambda.$$
This, along with the estimate (\ref{ineq:estimate-on-wc}) on $w_c$ and $\|w_o\|_0\leq C$, implies
$$\left|\int_{\mathbb T^d}|v_1|^2-|v|^2-|w_o|^2\,dx\right|\leq C\frac\mu{\lambda^{1-\alpha}}.$$
Now by definition (\ref{eq:rho(s)}) of $\rho$ we have
${\rm tr}\,\,R=d\rho=\frac1{(2\pi)^2}\left(e(t)\left(1-\frac\delta2\right)-\int_{\mathbb T^d}|v|^2\,dx\right)$.
Putting the above together finishes the proof.
\stopproof

\subsection{Estimates on the Reynolds stress}
In order to clarify the choice for $p_1$ as in (\ref{eq:p1}), we will temporarily write $p_1=p+q$.
With this, we have
\begin{eqnarray}
{\rm div}\,\mathring R_1&=&\partial_tv_1+{\rm div}\,(v_1\otimes v_1)+\nabla p_1\notag\\
&=&\partial_tw_o+v\cdot \nabla w_o\notag\\
&&+{\rm div}\,(w_o\otimes w_o+q{\rm Id}+\mathring R)\notag\\
&&+\partial_tw_c+{\rm div}\,(v_1\otimes w_c+w_c\otimes v_1-w_c\otimes w_c+v\otimes w_o).\notag
\end{eqnarray}
We split the Reynolds stress tensor into the {\it transport part}, the {\it oscillation part}, 
and the {\it error} as shown on the right-hand side of the above identity.
In the remainder of this Section we  estimate these terms separately.

\begin{lemma}[The transport part]
\begin{equation}
\|\mathcal R(\partial_tw_o+v\cdot\nabla w_o)\|_\alpha\leq C\left(\frac{\lambda^\alpha}\mu+\frac{\mu^2}{\lambda^{1-\alpha}}\right).
\end{equation}
\label{ineq:estimate-on-transport-part}
\end{lemma}
\startproof
This is proved in Lemma~7.1, p.~22 of \cite{DS3} for $d=3$ 
and the proof given there is valid for $d=2$ as well.
\stopproof

\begin{lemma}[The oscillation term]\label{lemma:estimate-on-oscillation-part}
\begin{equation}
\left\|\mathcal R\left({\rm div}\,(w_o\otimes w_o+q{\rm Id}+\mathring R)\right)\right\|_\alpha \leq C\frac{\mu^2}{\lambda^{1-\alpha}}.
\label{ineq:estimate-on-oscillation-part}
\end{equation}
\end{lemma}
\startproof
This is proved in Lemma~7.2, p.~23 of \cite{DS3} for $d=3$.
The main difference in the case $d=2$ is in the role of $q$.

Recalling that $R(y,s)=\rho(s){\rm Id}-\mathring R(y,s)$, see (\ref{eq:R(y,s)}),
and noting that $\rho=\rho(t)$ is a function of $t$ only,
\begin{eqnarray}
&&{\rm div}\,(w_o\otimes w_o+\mathring R+q{\rm Id}) \notag\\
=&&{\rm div}\,(w_o\otimes w_o-R+q{\rm Id})\notag\\
=&&{\rm div}\,\left(w_o\otimes w_o-R-\left(\frac{|w_o|^2}2-\nu\frac{\psi_o^2}2\right){\rm Id}\right)
+\nabla\left(q+\frac{|w_o|^2}2-\nu\frac{\psi_o^2}2\right)\notag\\
=&&{\rm div}\,_y\left(W\otimes W-R-\left(\frac{|W|^2}2-\nu\frac{\Psi^2}2\right){\rm Id}\right)\notag\\
&+&\lambda {\rm div}\,_\xi\left(W\otimes W-\left(\frac{|W|^2}2-\nu\frac{\Psi^2}2\right){\rm Id}\right)\notag\\
&+&\nabla\left(q+\frac{|w_o|^2}2-\nu\frac{\psi_o^2}2\right)\notag\\
=&&{\rm div}\,_y\left(W\otimes W-R-\left(\frac{|W|^2}2-\nu\frac{\Psi^2}2\right){\rm Id}\right)\notag\\
=&&{\rm div}\,_y\left(\sum_{|k|^2=\nu}U_ke^{i\lambda k\cdot x}-\sum_{|k|^2=\nu}\tilde a_ke^{i\lambda k\cdot x}{\rm Id}\right)\notag
\end{eqnarray}
where two cancelations occur by construction of $w_o$, see (\ref{eq:div(WoW)=nabla(P)}), 
and by definition of $q$, see (\ref{eq:p1}).
The estimate follows from Corollary~\ref{corollary:estimates-for-R} with $m=1$.
\stopproof

\begin{lemma}[The error - I]
\begin{equation}
\|\mathcal R\partial_tw_c\|_\alpha\leq C\frac{\mu^2}{\lambda^{1-\alpha}}.
\label{ineq:estimate-on-error-I}\end{equation}
\end{lemma}
\startproof
This is proved in Lemma~7.3, p.~23 of \cite{DS3} for $d=3$.
For clarity we briefly recall it for $d=2$ with the appropriate adjustments.

Recall that $u_c=\sum_{|k|^2=\nu} \psi_k(\lambda x)\nabla^\perp a_k(x,t,\lambda t)$ 
was defined in (\ref{eq:introduce-uc}) and thus
\begin{eqnarray}
\partial_t u_c(x,t)
=&&\lambda\sum_{|k|^2=\nu}\psi_k(\lambda x)\nabla^\perp\partial_\tau a_k(x,t,\lambda t)\notag\\
&+&\sum_{|k|^2=\nu}\psi_k(\lambda x)\nabla^\perp \partial_sa_k(x,t,\lambda t)\notag
\end{eqnarray}
But for any vector-valued function $A(y,s,\tau)$, we have
$$
{\rm div}\,\left(A(x,t,\lambda t)\otimes \frac k{|k|^3}e^{i\lambda k\cdot x}\right)
=i\lambda A(x,t,\lambda t)\frac{e^{i\lambda k\cdot x}}{|k|}+e^{i\lambda k\cdot x}\left(\frac k{|k|^3}\cdot\nabla\right)A(x,t,\lambda t)
$$
hence
$$\psi_k(\lambda x)A(x,t,\lambda t)e^{i\lambda k\cdot x}
=\frac1{i\lambda}{\rm div}\,\left(A(x,t,\lambda t)\otimes \frac k{|k|^3}e^{i\lambda k\cdot x}\right)
-\frac1{i\lambda}e^{i\lambda k\cdot x}\left(\frac k{|k|^3}\cdot \nabla\right) A(x,t,\lambda t).
$$
Therefore $\partial_tu_c$ is of the form
$$\partial_tu_c={\rm div}\,U_c+\tilde u_c$$
where
$$\|U_c\|_\alpha\leq C\mu\lambda^\alpha,\qquad \|\tilde u_c\|_\alpha\leq C\mu^2\lambda^\alpha$$
owing to the estimates from Lemma~\ref{lemma:estimates-on-coefficients}.
In turn,
\begin{eqnarray}
\|\mathcal R\partial_t w_c\|_\alpha
&\leq&\frac1\lambda\left(\|\mathcal R\mathcal Q{\rm div}\,U_c\|_\alpha+\|\mathcal R\mathcal Q\tilde u_c\|_\alpha\right)\notag\\
&\leq&\frac C\lambda\left(\|U_c\|_\alpha+\|\tilde u_c\|_\alpha\right)\notag\\
&\leq&C\frac{\mu^2}{\lambda^{1-\alpha}}\notag
\end{eqnarray}
owing to the estimates (\ref{ineq:Schauder-estimate-for-Q}), (\ref{ineq:Schauder-estimate-for-R}),  
and  (\ref{ineq:Schauder-estimate-for-RQdiv}).
\stopproof

\begin{lemma}[The error - II]
\begin{equation}
\|\mathcal R({\rm div}\,(v_1\otimes w_c+w_c\otimes v_1-w_c\otimes w_c))\|_\alpha\leq C\frac\mu{\lambda^{1-2\alpha}}.
\label{ineq:estimate-on-error-II}
\end{equation}
\end{lemma}
\startproof
This is proved in Lemma~7.4, p.~24 of \cite{DS3} for $d=3$
and the proof given there is valid in the case $d=2$ as well.
\stopproof

\begin{lemma}[The error - III]
\begin{equation}
\|\mathcal R({\rm div}\,(v\otimes w_o))\|_\alpha\leq C\frac{\mu^2}{\lambda^{1-\alpha}}.
\label{ineq:estimate-on-error-III}
\end{equation}
\end{lemma}
\startproof
This is proved in Lemma~7.4, p.~24 of \cite{DS3} for $d=3$
and we briefly indicate the proof in the case $d=2$.
Using that ${\rm div}\,_\xi b_k=0$, we find
\begin{eqnarray}
{\rm div}\,(v\otimes w_o)
&=&w_o\cdot\nabla v+({\rm div}\,w_o)v\notag\\
&=&\sum_{|k|^2=\nu}a_k(b_k\cdot \nabla)v+(\nabla a_k\cdot b_k)v\notag\\
&=&\sum_{|k|^2=\nu}\left[a_k\left(\frac{ik^\perp}{|k|}\cdot\nabla\right)v+\left(\nabla a_k\cdot \frac{ik^\perp}{|k|}\right)v\right]e^{i\lambda k\cdot x}.\notag
\end{eqnarray}
The estimate follows from Corollary~\ref{corollary:estimates-for-R} with $m=1$.
\stopproof

\section{Conclusions}
\subsection{Proof of Proposition~\ref{proposition:iteration}}
Recall that $e(t)$ is given as in Proposition~\ref{proposition:iteration}
and that $(v,p,\mathring R)$ is assumed to solve the Euler-Reynolds system (\ref{eq:Euler-Reynolds})
and to satisfy the bounds (\ref{ineq:bound-on-e(t)-int|v|2}) and (\ref{ineq:bound-on-Ro}).

We have now all estimates available in order to fix the parameters $\mu, \lambda$, and $\alpha$
so that the estimates (\ref{ineq:bound-on-e(t)-int|v1|2}), (\ref{ineq:bound-on-R1o}),
(\ref{ineq:bound-on-v1-v}), and (\ref{ineq:bound-on-p1-p}) may hold.
For simplicity we will take
\begin{equation}\mu=\lambda^\beta
\label{eq:mu=lambda^beta}
\end{equation}
for some $\beta$ to be determined
(although strictly speaking this can only hold up to some constant depending only on $\beta$ 
since it is required that $\mu\in\mathbb N$).
Recall that $C$ denotes a generic constant (possibly) depending on $e$, $v, \mathring R$, 
$\nu$, $\alpha$, and $\delta$, but not on $\lambda$ nor $\mu$.

Recall that the constant $M$ of Proposition~\ref{proposition:iteration} has already been fixed 
in (\ref{ineq:fix-parameter-M}) so that
\begin{equation}\max \big\{ \sqrt\nu\|\psi_o\|_0,\, \|w_o\|_0\big\}
\leq \frac{\sqrt M}2\sqrt\delta.
\label{ineq:fix-parameter-M-conclusion}
\end{equation}
Since $v_1-v=w_o+w_c$, the bound (\ref{ineq:bound-on-v1-v}) on $v_1-v$ follows provided
$$\|w_c\|_\alpha\leq C\frac \mu{\lambda^{1-\alpha}}=C\lambda^{\alpha+\beta-1}\leq \sqrt\delta.$$
The bound (\ref{ineq:bound-on-p1-p}) on $p_1-p=-\left(\frac{|w_o|^2}2-\nu\frac{\psi_o^2}2\right)$
also follows from (\ref{ineq:fix-parameter-M-conclusion}).

The bound (\ref{ineq:bound-on-e(t)-int|v1|2}) on the energy follows from (\ref{ineq:estimate-on-e(t)-int|v1|2})
provided 
$$C\frac\mu{\lambda^{1-\alpha}}=C\lambda^{\alpha+\beta-1}\leq \frac\delta8\min_{t\in[0,T]}e(t).$$

Finally, the estimates (\ref{ineq:estimate-on-transport-part}), (\ref{ineq:estimate-on-oscillation-part}),
(\ref{ineq:estimate-on-error-I}), (\ref{ineq:estimate-on-error-II}), and (\ref{ineq:estimate-on-error-III}), 
imply that 
$$\|\mathring R_1\|_\alpha\leq C(\lambda^{\alpha-\beta}+\lambda^{\alpha+2\beta-1}+\lambda^{2\alpha+\beta-1}).$$

In conclusion, imposing
\begin{equation}\alpha<\beta\qquad {\rm and}\qquad \alpha+2\beta<1
\label{ineq:conditions-on-alpha-and-beta}
\end{equation}
ensures that the bounds (\ref{ineq:bound-on-e(t)-int|v1|2}), (\ref{ineq:bound-on-R1o}),
(\ref{ineq:bound-on-v1-v}), and (\ref{ineq:bound-on-p1-p}) hold 
provided $\lambda$ is chosen sufficiently large.

This concludes the proof of Proposition~\ref{proposition:iteration}.

\subsection{A quantified version of Proposition \ref{proposition:iteration}}
In what follows we fix $d=3$.
The proof of Theorem~\ref{theorem:3D-flows-are-not-2D} follows from a closer
analysis of the proof of Theorem~\ref{theorem:dissipative-continuous}, in particular
from a more precise version of Proposition \ref{proposition:iteration} which will be stated below.
In order to prove it we collect first a series of estimates which make some of the statements
in the proof of Proposition \ref{proposition:iteration} more precise. 

\begin{lemma}\label{lemma:iteration}
Let $e, v, p, \mathring R$ be given.
Let $R, w_o, w_c, v_1, \mathring R_1$ be defined as in 
(\ref{eq:R(y,s)}), (\ref{eq:wo(x,t)}), (\ref{eq:v1=v+wo+wc}), and (\ref{eq:R1o}).
Then, there exist constants $C=C(e,v,p,\mathring R)$ and $\gamma>0$ such that
\begin{eqnarray}
\|w_c\|_\alpha&\leq& C\lambda^{-\gamma}\label{ineq:|wctilde|alpha}\\
\left|\int_{\mathbb T^3}v\cdot w_o\,dx\right|&\leq&C\lambda^{-\gamma}\label{ineq:int(vwotilde)}\\
\left|\int_{\mathbb T^3}w_o\otimes w_o\,dx-R\,dx\right|&\leq&C\lambda^{-\gamma}\label{ineq:int(wotildexwotilde-Rtilde)}\\
\left|\int_{\mathbb T^3}v_1\otimes v_1-v\otimes v-w_o\otimes w_o\,dx\right|&\leq&C\lambda^{-\gamma}\label{ineq:int(|vtilde|2-|v|2-|wotilde|2)}\\
\|\mathring R_1\|_\alpha&\leq&C\lambda^{-\gamma}.\label{ineq:|Rtildeo|alpha}
\end{eqnarray}
\end{lemma}
\startproof
Take $0<\alpha<\beta<1$ satisfying $\alpha+2\beta<1$, see (\ref{ineq:conditions-on-alpha-and-beta}), and set
$$\gamma=\min\left(1-\alpha-2\beta,\beta-\alpha,1-2\alpha-\beta\right).$$
Estimate (\ref{ineq:|wctilde|alpha}) is Lemma~\ref{lemma:estimate-on-the-corrector}.
Estimates (\ref{ineq:int(vwotilde)}), (\ref{ineq:int(wotildexwotilde-Rtilde)}) 
and (\ref{ineq:int(|vtilde|2-|v|2-|wotilde|2)}) follow from the proof of Lemma~\ref{lemma:estimate-on-the-energy}.
Specifically, (\ref{ineq:int(vwotilde)}) follows since $v$ is fixed and $w_o$ is oscillatory,
and (\ref{ineq:int(wotildexwotilde-Rtilde)}) follows since $w_o\otimes w_o- R$ is
the sum of oscillatory terms.
For (\ref{ineq:int(|vtilde|2-|v|2-|wotilde|2)}) we write
$$\int_{\mathbb T^3}v_1\otimes v_1\,dx
=\int_{\mathbb T^3}v\otimes v\,dx+\int_{\mathbb T^3}w_o\otimes w_o\,dx
+\int_{\mathbb T^3}\left(v\otimes w+w\otimes v
+w_o\otimes w_c+w_c\otimes w_o\right)\,dx.$$
The estimate follows since $w_o$ is oscillatory and $w_c$ is small ($v$ is fixed).
Estimate (\ref{ineq:|Rtildeo|alpha}) follows from Lemmas~7.1-7.4 of \cite{DS3}.
\stopproof

\begin{proposition}\label{proposition:quantified-iteration}
Let $e$ be as above and let $v, p, \mathring R$ solve the Euler-Reynolds system \eqref{eq:Euler-Reynolds}.
Suppose that there exist $0<\delta\leq 1$ satisfying
$$\left|e(t)(1-\delta)-\int_{\mathbb T^3}|v(x,t)|^2\,dx\right|\leq \frac\delta4 e(t)$$
$$\|\mathring R\|_0\leq \eta\delta.$$ %\min\left(\eta\delta,\varepsilon'\right).
Let $R$, $v_1=v+w=v+w_o+w_c$, $p_1=p+q$, $\mathring R_1$ be defined as in 
(\ref{eq:R(y,s)}), (\ref{eq:wo(x,t)}), (\ref{eq:v1=v+wo+wc}), (\ref{eq:p1}), and (\ref{eq:R1o}),
so that they satisfy the Euler-Reynolds system:
$$\partial_tv_1+{\rm div}\,(v_1\otimes v_1)+\nabla p_1={\rm div}\,\mathring R_1.$$
For any $\varepsilon'$, the following inequalities are satisfied provided $\lambda$ is chosen 
sufficiently large depending on $(v,p,\mathring R)$.
\begin{eqnarray}
\left|e(t)\left(1-\frac{\delta}2\right)-\int_{\mathbb T^3}|v_1(t)|^2\,dx\right|&\leq& \frac\delta8 e(t)
\label{ineq:improved-energy-bound}\\
\|\mathring R_1\|_0&\leq& \min\left(\frac12\eta\delta,\varepsilon'\right)
\label{ineq:bound-on-Rtildeo}\\
\|v_1-v\|_0&\leq& M\sqrt{\delta},\label{ineq:bound-on-|vtilde-v|0}\\
\sup_t\|v_1-v\|_{H^{-1}(\mathbb T^3)}&\leq &\varepsilon',\label{ineq:bound-on-|vtilde-v|H-1}\\
\|p_1-p\|_0&\leq&M\delta\label{ineq:bound-on-|ptilde-p|0}\\
\left|\int_{\mathbb T^3}v_1\otimes v_1\,dx
-\int_{\mathbb T^3}v\otimes v\,dx-\int_{\mathbb T^3}w_o\otimes w_o\,dx\right|
&\leq& \varepsilon'
\label{ineq:bound-integral-error-for-intvivj}\\
\left|\int_{\mathbb T^3} w_o\otimes w_o-R\,dx\right|&\leq& \varepsilon'.
\label{ineq:bound-on-int(woitildewojtilde-Rtildeij)}
\end{eqnarray}
\end{proposition}
\startproof
With
\begin{equation}\rho(t)=\frac1{3(2\pi)^3}\left[e(t)\left(1-\frac{\delta}2\right)-\int_{\mathbb T^3}|v(x,t)|^2\,dx\right]
\label{def:rhotilde}
\end{equation}
we have
\begin{equation}
\frac{\min_{0\leq t\leq T}e(t)\delta}{12(2\pi)^3}<\rho(t)\leq\frac{\max_{0\leq t\leq T}e(t)\delta}{4(2\pi)^3}.
\label{ineq:rhotilde(t)}
\end{equation}
Owing to the definition (\ref{def:eta}) of $\eta$ we obtain
$$\left\|\frac{R}{\rho}-{\rm Id}\right\|_0
=\left\|\frac{\mathring R}{\rho}\right\|_0
< r_0.
$$
Thus, $v_1, p_1, R_1$ etc. can be defined and estimated as in Lemma~\ref{lemma:iteration}.
We will now choose $\lambda$ in Lemma~\ref{lemma:iteration} sufficiently large (depending on $e, v, p, \tilde R$) 
so that the desired estimates hold.
\medskip

Recalling the definition (\ref{eq:R(y,s)}) 
of $R$ in terms of $\rho$, we have
\begin{eqnarray}
&&e(t)\left(1-\frac{\delta}2\right)-\int_{\mathbb T^3}|v_1(t)|^2\,dx\notag\\
&=&e(t)\left(1-\frac{\delta}2\right)-\int_{\mathbb T^3}|v(t)|^2\,dx
-\int_{\mathbb T^3} |w_o|^2\,dx + \int_{\mathbb T^3}\left(|v_1|^2-|v|^2-|w_o|^2\right)\,dx\notag\\
&=&-\int_{\mathbb T^3}\left(|w_o|^2 - {\rm tr}\,R\right)\,dx
+\int_{\mathbb T^3}\left(|v_1|^2-|v|^2-|w_o|^2\right)\,dx.\notag
\end{eqnarray}
Invoking (\ref{ineq:int(wotildexwotilde-Rtilde)}) and (\ref{ineq:int(|vtilde|2-|v|2-|wotilde|2)}),
we obtain (\ref{ineq:improved-energy-bound}) with large $\lambda$.

From (\ref{ineq:|Rtildeo|alpha}) it is clear that the bound (\ref{ineq:bound-on-Rtildeo}) on $\mathring R_1 $ 
holds with $\lambda$ large.

From (\ref{ineq:|wctilde|alpha}), $w_c$ can be made arbitrarily small with $\lambda$ large,
while from (\ref{eq:wo(x,t)}) 
and (\ref{def:M'}) we have $|w_o|\leq \frac12M'\sqrt{\tilde\rho}\leq \frac12M\sqrt{\delta}$.
Thus, the bounds (\ref{ineq:bound-on-|vtilde-v|0}) on $v_1-v=w_o+w_c$
and (\ref{ineq:bound-on-|ptilde-p|0}) on $p_1-p=-\frac{|w_o|^2}2$ follow with $\lambda$ large.
Furthermore, $w_o$ is a fast oscillating function and $w_c$ is small in the sup-norm,
so that we achieve the $H^{-1}$-bound (\ref{ineq:bound-on-|vtilde-v|H-1}) on $v_1-v$ as well with $\lambda$ large.

Finally, (\ref{ineq:int(wotildexwotilde-Rtilde)}) and (\ref{ineq:int(|vtilde|2-|v|2-|wotilde|2)})
trivially imply (\ref{ineq:bound-integral-error-for-intvivj}) and (\ref{ineq:bound-on-int(woitildewojtilde-Rtildeij)})
with large $\lambda$.
\stopproof

\subsection{Proof of Theorem~\ref{theorem:3D-flows-are-not-2D}}
We are now ready to prove that the solution $v,p$ constructed in Theorem~\ref{theorem:dissipative-continuous}
in the case $d=3$ satisfies the estimates of Theorem~\ref{theorem:3D-flows-are-not-2D}
provided $\lambda$ is chosen sufficiently large at each iteration.
Set
$$\varepsilon_0:=\frac12\,\min_{0\leq t\leq T}e(t).$$
By choosing $\varepsilon>0$ smaller, if necessary, we may assume that $\varepsilon<\frac{1}{6}\varepsilon_0$. 
We recall that the solution $(v,p)$ is obtained as a limit
$$v:=\lim_n v^{(n)},\qquad p:=\lim_np^{(n)}$$
where the sequences $v^{(n)}$, $p^{(n)}$ are as follows.
We let 
$$v^{(0)}=0, \qquad p^{(0)}=0, \qquad \mathring R^{(0)}=0.$$
For $n\geq 0$, we construct 
$$v^{(n+1)}=v^{(n)}+w^{(n+1)},\qquad p^{(n+1)}=p^{(n)}+q^{(n+1)},\qquad \mathring R^{(n+1)}$$
using Proposition~\ref{proposition:quantified-iteration} 
with $v^{(n)}=v$, $p^{(n)}=p$, $\mathring R^{(n+1)}=\mathring R_1$, $v^{(n+1)}=v_1$, $p^{(n+1)}=p_1$
\begin{equation}
\delta=2^{-n}, %\qquad \delta=\delta^{(0)}2^{-n-1}
\qquad\mbox{and}\qquad \varepsilon'=\varepsilon 2^{-n-1},
\label{eq:use-Prop-for-iteration-n-to-n+1}
\end{equation}
so that
$$
\partial_t v^{(n)}+{\rm div}\,(v^{(n)}\otimes v^{(n)})+\nabla p^{(n)}={\rm div}\,\mathring{R}^{(n)}.
$$
Recall that, in Proposition~\ref{proposition:quantified-iteration}, the velocity $v_1$ is obtained
by adding two perturbation $w_o$ and $w_c$ to $v$, defined, respectively, in \eqref{eq:wo(x,t)}
and in \eqref{eq:v1=v+wo+wc}. We use the notation $w_o^{(n+1)}=w_o$ and $w_c^{(n+1)}=w_c$.
Moreover, we set 
\[
\rho^{(n)} (t) = \frac1{3(2\pi)^3}\left(e(t)\left(1-2^{-n-1}\right)-\int_{\mathbb T^3}|v^{(n)}|^2(x,t)\,dx\right)
\]
and $R^{(n)}=\rho^{(n)} {\rm Id} - \mathring{R}^{(n)}$ (observe that therefore 
$\rho^{(n)}=\rho$ and $\mathring{R}^{(n)}=R$, with $\rho$ and $R$ as in
\eqref{eq:rho(s)} and \eqref{eq:R(y,s)}).
In particular,
$$\rho^{(0)}(t)
\geq \frac{\varepsilon_0}{3(2\pi)^3},\qquad t\in [0,T]$$
according to (\ref{def:rhotilde}) with $v=0$ and $\delta=1$ (step $n=0$).

Now
$$\|v\|_{H^{-1}}\leq \sum_{n=0}^\infty\|w^{(n+1)}\|_{H^{-1}}\leq \varepsilon$$
easily follows from (\ref{ineq:bound-on-|vtilde-v|H-1}).
This establishes (\ref{ineq:|v|H-1<epsilon}).

Next, for each $n$,
\begin{eqnarray}
\int_{\mathbb T^3} v^{(n+1)}\otimes v^{(n+1)}\,dx
&=&\int_{\mathbb T^3} v^{(n)}\otimes v^{(n)} +\int_{\mathbb T^3} R^{(n)}\,dx\notag\\
&&+\int_{\mathbb T^3}\left(w_o^{(n)}\otimes w_o^{(n+1)}-R^{(n)}\right)\,dx\notag\\
&&+\int_{\mathbb T^3} \left(v^{(n+1)}\otimes v^{(n+1)} - v^{(n)}\otimes v^{(n)} - w_o^{(n+1)}\otimes w_o^{(n+1)}\right)\,dx\notag\\
&=&\int_{\mathbb T^3} v^{(n)}\otimes v^{(n)} + (2\pi)^3\rho^{(n)}{\rm Id}-\int_{\mathbb T^3}\mathring R^{(n)}\,dx\notag\\
&&+\int_{\mathbb T^3}\left(w_o^{(n+1)}\otimes w_o^{(n+1)}-R^{(n)}\right)\,dx\notag\\
&&+\int_{\mathbb T^3} \left(v^{(n+1)}\otimes v^{(n+1)} - v^{(n)}\otimes v^{(n)} - w_o^{(n+1)}\otimes w_o^{(n+1)}\right)\,dx.\notag
\end{eqnarray}
But the sum of the last two integrals is bounded by $\varepsilon 2^{-n}$ 
using (\ref{ineq:bound-integral-error-for-intvivj}), (\ref{ineq:bound-on-int(woitildewojtilde-Rtildeij)}),
and (\ref{eq:use-Prop-for-iteration-n-to-n+1}). Therefore, taking the trace we obtain
\begin{align*}
\int_{\mathbb T^3}|v^{(n+1)}(x,t)|^2\,dx&=\int_{\mathbb T^3}|v^{(n)}(x,t)|^2\,dx+3(2\pi)^3\rho^{(n)}(t)+O(\varepsilon 2^{-n})\\
&=e(t)(1-2^{-n-1})+O(\varepsilon 2^{-n}),
\end{align*}
where $O(\varepsilon 2^{-n})$ denotes an error term bounded by $\varepsilon 2^{-n}$.
In particular this yields $\rho^{(n)}(t)=\frac{1}{3(2\pi)^3}2^{-n-1}e(t)+O(\varepsilon 2^{-n})$.
Furthermore, by summing over $n$ we obtain
\begin{eqnarray}
\int_{\mathbb T^3} v\otimes v\,dx
&=&-\sum_{n=0}^{\infty}\int_{\mathbb T^3} \mathring R^{(n)} \,dx
+ (2\pi)^3\sum_{n=0}^{\infty}\rho^{(n)}{\rm Id} + O(\varepsilon).\notag
\end{eqnarray}
Since, by (\ref{ineq:bound-on-Rtildeo}) and (\ref{eq:use-Prop-for-iteration-n-to-n+1}) we also have
$$\left|\sum_{n=0}^{\infty}\int_{\mathbb T^3} \mathring R^{(n)}\,dx \right|\leq \varepsilon$$
we deduce
$$
\int_{\mathbb T^3} v\otimes v\,dx =\frac{1}{3}e(t){\rm Id} + O(\varepsilon)\,.
$$
\stopproof

%%%%%%%%%%%%%%%%%%%%%%%%%%%%%%%%%%%%

\bibliographystyle{amsalpha}
\bibliography{eulerbib}

\end{document}